\documentclass[11pt]{article}
\usepackage[utf8]{inputenc}
\usepackage[a4paper]{geometry}
\usepackage[english]{babel}
\usepackage[super]{nth}
\usepackage{amsmath, amsfonts, amsthm, physics, enumitem}
\usepackage{titlesec}
\usepackage{tikz} 
\usetikzlibrary{positioning}
\usepackage{mathtools}

\usetikzlibrary{patterns,arrows,decorations.markings}
\usetikzlibrary{shapes.multipart,calc,shapes.geometric,positioning}
\usepackage{authblk}
\usepackage{tikz-cd}
\usepackage{centernot}

\usepackage{tikz}
\usetikzlibrary{positioning}
\usepackage{centernot}
\usepackage{caption}
\usepackage{float}
\usepackage{CJKutf8}
\usepackage[rightcaption]{sidecap}
\usepackage{mathtools}

\pgfdeclareimage[height=15pt]{g_tick}{green_tick}
\pgfdeclareimage[height=15pt]{y_tick}{yellow_tick}
\pgfdeclareimage[height=15pt]{r_cross}{Red_X}

\newtheorem{thm}{Theorem}[section]

\newtheorem{q}[thm]{Question}

\theoremstyle{definition}

\def \p{\partial}

\newcommand{\midarrow}{\tikz \draw[-triangle 90] (0,0) -- +(.1,0);}

\newcommand{\KT}{\mathrm{KT}^{4}}

\newcommand{\der}[1]{\frac{\partial}{\partial#1}}
\newcommand{\dif}[1]{\mathrm{d}#1}

\makeatletter
\newcommand{\subjclass}[2][2020]{%
  \let\@oldtitle\@title%
  \gdef\@title{\@oldtitle\footnotetext{#1 \emph{Mathematics subject classification.} #2}}%
}
\newcommand{\keywords}[1]{%
  \let\@@oldtitle\@title%
  \gdef\@title{\@@oldtitle\footnotetext{\emph{Key words and phrases.} #1.}}%
}
\makeatother

\title{Almost complex Hodge theory}
\author{ Weiyi Zhang\thanks{Weiyi.Zhang@warwick.ac.uk}}
\affil{Mathematics Institute,  University of Warwick, Coventry, CV4 7AL, England}

\date{}

\begin{document}



\maketitle

\begin{abstract}
We review the recent development of Hodge theory for almost complex manifolds. This includes the determination of whether the Hodge numbers defined by $\bar\partial$-Laplacian are almost complex, almost K\"ahler, or birational invariants in dimension four. 

\end{abstract}

\tableofcontents

\section{Introduction}
Hodge theory is a method introduced by Hodge in the 1930s to study the cohomology groups of compact manifolds using the theory of elliptic partial differential equations. Not only has Hodge theory since become part of the standard repertoire in algebraic geometry, particularly through its connection to the study of algebraic cycles, but also the elliptic theory has become a fundamental tool to study the topology of manifolds. 

The most fundamental idea in classical Hodge theory for complex manifolds is the introduction of the finite dimensional vector spaces of $\bar\partial$-harmonic $(p, q)$-forms $\mathcal H^{p, q}$ with respect to a K\"ahler (or Hermitian) metric. Each space $\mathcal H^{p, q}$ can be identified with a coherent sheaf cohomology group, called the Dolbeault group, which depends only on the underlying complex manifold but not on the choice of the Hermitian metric. Their dimensions $h^{p, q}$, called the Hodge numbers, are important invariants of complex manifolds.  They do not change when the complex structures are K\"ahler and vary continuously, but they are in general not topological invariants. When the manifold is K\"ahler, these spaces give rise to a decomposition of the cohomology groups with complex coefficients. 

As discussed in \cite{Hir} and recalled in Section \ref{hodgeno}, we can still define the group $\mathcal H^{p, q}$ for closed almost complex manifolds with Hermitian metrics.   
Kodaira and Spencer asked whether this almost complex version of $\mathcal H^{p, q}$ would depend on the choice of Hermitian metrics, which appeared as Problem 20 in Hirzebruch's 1954 problem list \cite{Hir}. 
 
\begin{q}[Kodaira-Spencer]\label{KS}
Let $M$ be a compact almost complex manifold. For any given Hermitian structure we can consider the numbers $h^{p,q}$. Are these numbers independent of the choice of the Hermitian structure?
\end{q} 

This is a question guiding the study of the almost complex Hodge theory. If for some pair $(p,q)$ the answer to this question is positive, then we have almost complex invariant $h^{p,q}$. In some cases, we can also show it is a birational invariant.  On the other hand, if the question is answered negatively for some $(p,q)$, then it would bring us some phenomena which are very different from classical Hodge theory. 

Although the almost complex Hodge theory could be very useful, in particular for geometrically interesting almost complex structures, not much was known beyond the attempts to develop harmonic theory for almost K\"ahler manifolds by Donaldson \cite{Don}, for strictly nearly K\"ahler $6$-manifolds by Verbitsky \cite{Ver}, and very recently the introduction of a variant of $\mathcal H^{p,q}$ using $\bar\partial$-$\mu$-harmonic forms by Cirici and Wilson \cite{CW}. Moreover, few non-trivial examples of $\mathcal H^{p,q}$ for non-integrable almost complex structures have been computed.

In this article, we review recent development of almost complex Hodge theory, which in particular includes fairly thorough solution for  Question \ref{KS} in dimension $4$ for all pairs $(p,q)$. Results on two related versions of harmonic $(p,q)$-forms, the Bott-Chern and $d$-harmonic forms, are also reviewed. These are done in Section \ref{hodgeno}. In Section \ref{KTex}, explicit calculations of Hodge numbers for a family of almost complex structures on the Kodaira-Thurston manifold are summarized. In Section \ref{pdeodent}, the from PDE to ODE and Number Theory method for the computation of $h^{0,1}$, which in particular leads to the negative answer to  Question \ref{KS} in \cite{HZ, HZ2}, is carefully explained.

This paper is based on a talk delivered at the conference {\it Cohomology of complex manifolds and special structures, II} in July 2021, although some results in \cite{Ho} after the conference are also contained.

\section{Hodge numbers on almost complex manifolds}\label{hodgeno}
We summarize the results for Hodge numbers on almost complex manifolds. 

\subsection{Preliminaries}
The almost complex structure $J$ induces the decomposition $$T^*M\otimes_{\mathbb R} \mathbb C=(T^*M)^{1, 0}\oplus (T^*M)^{0,1}$$
and in turn $$\Lambda^rT^*M\otimes_{\mathbb R} \mathbb C=\oplus_{p+q=r}\Lambda^{p, q}M.$$

As in integrable complex case, we define
$$\bar{\partial}=\pi^{p,q+1}\circ d: \Omega^{p,q}(M)\rightarrow \Omega^{p,q+1}(M)$$
$$\partial=\pi^{p+1,q}\circ d: \Omega^{p,q}(M)\rightarrow \Omega^{p+1,q}(M),$$
where $\pi^{p,q}$ be the projection to $\Lambda^{p, q}M$ and $\Omega^{p, q}(M)=\Gamma(M, \Lambda^{p,q})$.

In general, we have $$d=\bar\mu+\bar\partial +\partial+\mu,$$ where $\mu$ and $\bar\mu$ have bidegrees $(2, -1)$ and $(-1,2)$ respectively. The almost complex structure is integrable if and only if $d=\bar\partial +\partial$.

A manifold is almost Hermitian provided it is almost complex and has a Riemannian metric $g$ compatible with $J$, {\it i.e.} $g(Ju, Jv)=g(u,v)$, $\forall u, v \in TM$. It is equivalent to an Hermitian metric $h$ on $TM^{1,0}\subset TM\otimes \mathbb C$ which is related to the Riemannian metric $g$ by $h=g-i\omega$, where $\omega(u, v)=g(Ju, v)$. 

The $*$ operator on an almost Hermitian manifold is the unique $\mathbb{C}$-linear operator $$*: \Lambda^{p,q}\rightarrow \Lambda^{n-q,n-p}$$ satisfying \begin{align*} g(\varphi_1, \varphi_2)dV=\varphi_1\wedge\overline{*\varphi_2}\end{align*} where $dV$ is the volume form of $g$ and $\varphi_1, \varphi_2\in \Lambda^{p,q}.$

In this paper, we assume $M$ is compact. Define an inner product on $\Omega^{p,q}(M)$ by $\langle \varphi_1, \varphi_2\rangle=\int_M g(\varphi_1,\varphi_2)dV$. Then $$\bar{\partial}^*=-*\partial*$$ is the formal $L^2$-adjoint of $\bar\partial$, {\it i.e.}  $$\langle \bar{\partial}\varphi_1, \varphi_2\rangle=\langle \varphi_1, \bar{\partial}^*\varphi_2\rangle.$$ 

The $\bar{\partial}$-Laplacian
$$\Delta_{\bar\partial}=\bar\partial\bar\partial^*+\bar\partial^*\bar\partial$$ is an elliptic operator. 

The space $\mathcal H^{p,q}\subset \Omega^{p,q}$ is defined to be the kernel of $\Delta_{\bar\partial}$. We have
  $$\mathcal H^{p,q}=\ker\Delta_{\bar\partial}=\ker\bar\partial\cap\ker\bar\partial^*= \ker\bar\partial\cap \ker\partial *.$$
 This would depend on the choice of the almost Hermitian metric, but we will omit the appearance of $h$ or $g$ for most part of this article except in \eqref{20b}. 
This setup could be generalized to Hermitian bundles with a pseudoholomorphic structure, {\it c.f.} \cite{CZ}. 

The focus of this paper is on Hodge numbers $h^{p,q}:=\dim\mathcal H^{p,q}$. Question \ref{KS} asks whether $h^{p,q}$ is independent of the choice of the almost Hermitian metric. If it is true for some pair $(p,q)$, $h^{p,q}$ defines an almost complex invariant. 

There are two related versions of harmonic $(p,q)$-forms. The first is $\mathcal H_d^{p,q}=\ker \Delta_d|_{\Omega^{p,q}}$, where $\Delta_d=dd^*+d^*d$. 
The second is a corresponding Bott-Chern version, defined by Piovani-Tomassini \cite{PT} based on the work of Kodaira and Spencer. It is $\mathcal H_{BC}^{p,q}=\ker \Delta_{BC}|_{\Omega^{p,q}}$, where $\Delta_{BC}=\partial\bar\partial\bar\partial^*\partial^*+\bar\partial^*\partial^*\partial\bar\partial+\partial^*\bar\partial\bar\partial^*\partial+\bar\partial^*\partial\partial^*\bar\partial+\partial^*\partial+\bar\partial^*\bar\partial$. We denote their dimensions by $h_d^{p,q}$ and $h_{BC}^{p,q}$ respectively. We will also mention relevant results on these two spaces in this article. 

We begin with results where $h^{p,q}$ behaves like in the complex setting. 

\subsection{Serre duality}
Serre duality still holds on compact almost complex manifolds in the sense that $$\mathcal{H}^{p,q} = *\overline{\mathcal{H}^{n-p,n-q}}$$ for any given almost Hermitian metric. Therefore, to know the whole Hodge diamond for an almost Hermitian $4$-manifold, we only need to compute $h^{1,0}$, $h^{2, 0}$, $h^{0,1}$, and $h^{1,1}$.

\subsection{$\boldsymbol{h^{p,0}}$ is an almost complex and birational invariant}

As $\bar\partial^*=0$ on $\Omega^{p,0}(M)$, $$\mathcal H^{p,0}= \ker \Delta_{\bar\partial}=\ker \bar\partial.$$
So $h^{p,0}$ is independent of the choice of the Hermitian structure. Apparently, $h^{0,0}=1$ when $M$ is compact and connected.

As a side remark, the identification also works for the bundle valued forms. In particular, the pseudoholomorphic plurigenera $$\dim H^0(M,\mathcal K^{\otimes m})=\dim \{s\in \Gamma(M, \mathcal K^{\otimes m}):  \bar\partial_ms=0\}, \quad \mathcal K=\Lambda^{n,0}$$ are finite numbers  and thus almost complex invariants, which is the starting point of my work with Haojie Chen on almost complex Kodaira dimension \cite{CZ, CZ2}.

In algebraic geometry, $h^{p,0}$ alongside plurigenera and Kodaira dimension are very important birational invariants. They can be extended to almost complex setting. First we need to introduce birational morphisms in the almost complex setting as follows. We define two almost complex manifolds $M$ and $N$ to be birational to each other if there are almost complex manifolds $M_1, \cdots, M_{n+1}$ and $X_1, \cdots, X_{n}$ such that $M_1=M$ and $M_{n+1}=N$, and there are degree one pseudoholomorphic maps $f_i: X_i\rightarrow M_i$ and $g_i: X_{i}\rightarrow M_{i+1}$, $i=1, \cdots, n$.

This extension is naturally suggested by the following theorem in \cite{ZCJM}.

\begin{thm}
Let $u: (X, J) \rightarrow (M, J_M)$ be a degree one pseudoholomorphic map between closed almost complex $4$-manifolds such that $J$ is almost K\"ahler. Then there exists a subset $M_1\subset M$, consisting of finitely many points, with the following significance:
\begin{enumerate}
\item The restriction $u|_{X\setminus u^{-1}(M_1)}$ is a diffeomorphism. 
\item At each point of $M_1$, the preimage is an exceptional curve of the first kind.  
\item $X\cong M\#k\overline{\mathbb CP^2}$ diffeomorphically, where $k$ is the number of irreducible components of the $J$-holomorphic $1$-subvariety $u^{-1}(M_1)$.
\end{enumerate}
\end{thm}

Here, an exceptional curve of the first kind is a connected $J$-holomorphic $1$-subvariety whose configuration is equivalent to the empty set through topological blowdowns.

Under this notion of birationality, it is shown in \cite{CZ} that $h^{p,0}$, $0\le p\le 2$, are birational invariants for closed almost complex $4$-manifolds by establishing the Hartogs's extension theorem for almost complex manifolds. 

We end this subsection with two remarks on almost complex birational invariants. First, the proof of birational invariance of $h^{p,0}$ is extendable to higher dimensions once we can show that a closed subset in a closed $2k$-dimensional almost complex manifold with finite $(2k-2)$-dimensional Hausdorff measure support a $J$-holomorphic subvariety
if it is endowed with ``positive cohomology assignment", which is a notion of intersection number with each local open disk.  
When $k=4$, it is a result of Taubes \cite{T}. Second, the method could be used to show that the pseudoholomorphic plurigenera, almost complex Kodaira dimension \cite{CZ}, and the dimension of  closed $J$-anti-invariant $2$-forms on a $4$-manifold \cite{BZ} are birational invariants. Recently, Holt \cite{Ho} shows that $h^{p,0}_{BC}$ are also birational invariants for closed almost complex $4$-manifolds.

A precise computation of $h^{1,0}$ and $h^{2,0}$ for a family of almost complex manifolds on the Kodaira-Thurston manifold will be provided in Section \ref{KTex}.

\subsection{$\boldsymbol{h^{0,1}}$ is not an almost complex nor an almost K\"ahler invariant}
However, Hodge numbers do not always behave as in the integrable case. In dimension $4$, the Hodge number $h^{0,1}$ behaves most drastically different from the integrable case. First, $h^{0,1}$ is no longer metric independent, which answers Question \ref{KS} negatively.

\begin{thm}[Holt-Z. \cite{HZ, HZ2}]\label{intro3}
There exist almost complex structures on the Kodaira-Thurston manifold such that $h^{0,1}$ varies with different choices of almost K\"ahler metrics (or non almost K\"ahler Hermitian metrics).
\end{thm}

We can also see how $h^{0,1}$ would change when we vary $J$. First, the classical Hodge numbers are constant in a small neighbourhood of the moduli of a given K\"ahler manifold \cite{Voi}, in particular they do not change when the complex structures are K\"ahler and vary continuously. Second, even if Hodge numbers are in general not topological invariants, we know they are bounded, for example by the Betti numbers, for a fixed compact complex manifold with K\"ahler structures. 

However, neither of the above statements for classical Hodge numbers is true in the almost complex setting when we vary almost complex structures in an almost K\"ahler family. 

\begin{thm}[Holt-Z. \cite{HZ}]\label{intro1}
There is a continuous family of non-integrable almost complex structures $J_b$, $b\in \mathbb R\setminus\{0\}$, on the Kodaira-Thurston manifold whose $h^{0,1}_{J_b}=h^{2,1}_{J_b}$ are computed using certain almost K\"ahler metrics. Then for any $n\in\mathbb Z^+$ such that $8\centernot\mid n$, there is a $b$ such that $h^{0,1}_{J_b}=n$.
\end{thm}

We will sketch the proofs of the above results in Section \ref{pdeodent}.

Notice these examples are all on the Kodaira-Thurston manifold. It would be interesting to know almost Hermitian metrics with different $h^{0,1}$ on other manifolds. These examples are expected on other nil or solv manifolds. What about simply connected $4$-manifolds? 

In higher dimensions, we choose $M=\KT\times \mathbb CP^n$ with the product almost complex structure $J_b\times J_{std}$ and product almost K\"ahler form $\omega_b\times \Omega_{FS}$.  We then have exactly the same statements of Theorem \ref{intro1} for this setting.

\subsection{$\boldsymbol{h^{1,1}}$ and its variants}
If we choose the Hermitian metric $h=g-i\omega$ to be almost K\"ahler, {\it i.e.} $d\omega=0$, we have the following generalized Hodge Index Theorem \cite{HZ}.
\begin{thm}[Holt-Z.]\label{akh11}
For any closed almost K\"ahler $4$-manifold $(M, J)$, $h^{1,1}$ is independent of the choice of almost K\"ahler metrics compatible with $J$. More precisely, $h^{1,1}=b^-+1$.
\end{thm}

Recently, Tardini and Tomassini \cite{TT} show that $h^{1,1}$ is not an almost complex invariant though. 
More precisely, they show that for any closed almost complex $4$-manifold $(M, J)$, 
\begin{equation}\label{tt}
h^{1,1}=\left\{\begin{array}{cl}
b^-+1 &\hbox{ globally conformal almost K\"ahler metric},\\
b^-& \hbox{ strictly locally conformal almost K\"ahler metric}.\\
\end{array}\right.
\end{equation}

The first case is a slight extension of Theorem \ref{akh11} using the fact that $h^{p,n-p}$ is a conformal invariant of Hermitian metrics on almost complex manifolds of dimension $2n$. In the second case, a strictly locally conformal metric is given by a non-degenerate $(1,1)$-form $\omega$ such that $d\omega=\alpha\wedge\omega$ where $\alpha$ is a $d$-closed, non $d$-exact, $1$-form. 

Recently, Tom Holt  shows that $b^-$ and $b^-+1$ are the only values which could be achieved by $h^{1,1}$ on a closed almost Hermitian $4$-manifold \cite{Ho}. On the other hand, for Bott-Chern cohomology, $h^{1,1}_{BC}$ is always $b^-+1$ for any closed almost Hermitian $4$-manifold \cite{Ho, PT}.

When $(M, J)$ is a complex surface, $h^{1,1}$ is independent of Hermitian metrics. Moreover, $h^{1,1}=b^-+1$ if and only if $(M, J)$ is K\"ahler. This is the well-known K\"ahler criterion for complex surfaces that a compact complex surface is K\"ahler if and only if the first Betti number $b_1$ is even, which is also equivalent to  $h^{1,1}=b^-+1$. This follows from Kodaira's classification of compact complex surfaces and the work of Miyaoka \cite{Miy} and Siu \cite{Siu}. There are also classification-free analytic proofs by Buchdahl \cite{Bu} and Lamari \cite{La}.

For an almost complex $4$-manifold, we hope to find similar almost K\"ahler criteria, although the one using the parity of $b_1$ apparently no longer works. However, it is reasonable to expect to distinguish almost K\"ahler ones from almost complex structures using some version of $(1,1)$ cohomology groups. 

Such a criterion exists  for the dimension of the space of $d$-harmonic $(1,1)$-forms $h_d^{1,1}$. Holt shows that on a closed almost Hermitian $4$-manifold $(M, J, \omega)$, $h_d^{1,1}=b^-+1$ if and only if $\omega$ is in the conformal class of an almost K\"ahler metric (and otherwise $h_d^{1,1}=b^-$). On the other hand, the corresponding statement for $h^{1,1}$ no longer holds. Piovani and Tomassini \cite{PT2} recently provide examples of closed almost Hermitian $4$-manifolds $(M, J, \omega)$ with $h^{1,1}=b^-+1$ but $\omega$ is not locally conformally almost K\"ahler. However, as pointed out to the author by Tom Holt, their almost complex structures $J$ are all almost K\"ahler. 

Hence, we might still hope to distinguish almost K\"ahler $J$ by looking at all compatible almost Hermitian structures. For example, we define an invariant for $(M, J)$  
\begin{equation}\label{20b}
\tilde h^{p,q}=\max_{(g, J)\hbox{\footnotesize{ almost Hermitian}}}h^{p,q}_{\bar\partial, g}.
\end{equation}
It is not yet known whether $\tilde h^{p,q}$ is always finite in general. For $4$-manifolds, the above mentioned result of Holt shows that $\tilde h^{1,1}=b^-$ or $b^-+1$. It is natural to ask whether $\tilde h^{1,1}$ characterizes almost K\"ahler structures, {\it i.e.},
 is it true that $\tilde h^{1,1}=b^-+1$ if and only if $J$ is almost K\"ahler? Recently, Piovani \cite{P} shows that this is not true by a nice computation on compact quotients of Lie groups.

\subsection{Hodge diamond}
To summarize Section $2$, we have the following diagram for the Hodge diamond of almost complex $4$-manifolds.

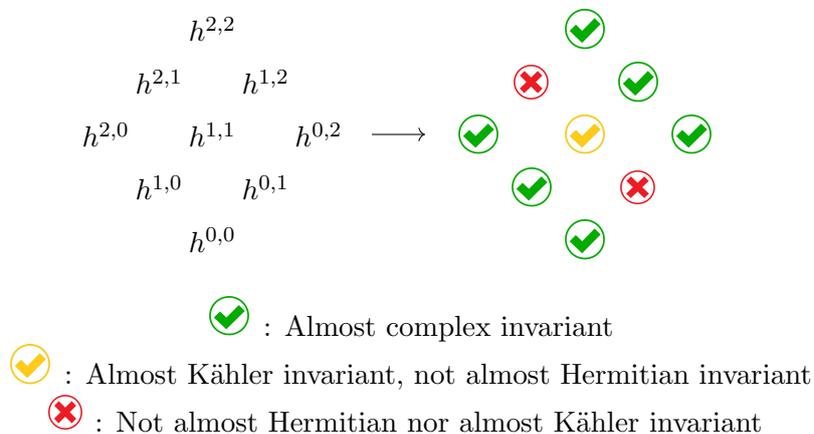
\begin{center}
\begin{figure}[H]
\centering
\captionsetup{justification=centering}
    \tikz [scale = .7,decoration={
       markings,
       mark=at position 1 with {\arrow[scale=3,black]{latex}};}] {
    
     \node (n1) at (0,0)  {$h^{0,0}$} ;
     \node (n2) at (1,1) {$h^{0,1}$};
     \node (n3) at (-1,1)  {$h^{1,0}$};
     \node (n4) at (0,2)  {$h^{1,1}$};
     \node (n5) at (2,2)  {$h^{0,2}$} ;
     \node (n6) at (-2,2) {$h^{2,0}$};
     \node (n7) at (-1,3)  {$h^{2,1}$};
     \node (n8) at (1,3)  {$h^{1,2}$};
     \node (n8) at (0,4)  {$h^{2,2}$};

    \draw [->] (3,2) -- (4,2);

    \node (n1) at (7,0)  {\pgfuseimage{g_tick}} ;
    \node (n2) at (8,1) {\pgfuseimage{r_cross}};
    \node (n3) at (6,1)  {\pgfuseimage{g_tick}};
    \node (n4) at (7,2)  {\pgfuseimage{y_tick}};
    \node (n5) at (9,2)  {\pgfuseimage{g_tick}} ;
    \node (n6) at (5,2) {\pgfuseimage{g_tick}};
    \node (n7) at (6,3)  {\pgfuseimage{r_cross}};
    \node (n8) at (8,3)  {\pgfuseimage{g_tick}};
    \node (n8) at (7,4)  {\pgfuseimage{g_tick}};
}
\caption*{   {\pgfuseimage{g_tick}}: Almost complex  invariant \\  {\pgfuseimage{y_tick}}: Almost K\"ahler  invariant, not  almost Hermitian  invariant \\  {\pgfuseimage{r_cross}} : Not almost Hermitian  nor almost K\"ahler invariant}
\end{figure}
\end{center}

By Serre duality and $h^{1,1}=b^-$ or $b^-+1$, $\tilde h^{0,1}=\tilde h^{2,1}$ is the only $\tilde h^{p,q}$ which is not known to be finite yet in for an almost complex $4$-manifold $(M, J)$. 

Holt extends the computation of $h^{0,1}$ to Bott-Chern cohomology for $h^{2,1}_{BC}$ and $h^{1,2}_{BC}$. He shows that they also vary with almost Hermitian metrics \cite{Ho}.  Combining discussions in this section and work of \cite{PT}, we have

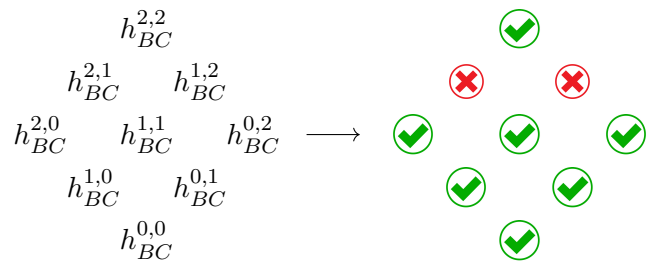
\begin{center}
\begin{figure}[H]
\centering
\captionsetup{justification=centering}
    \tikz [scale = .7,decoration={
       markings,
       mark=at position 1 with {\arrow[scale=3,black]{latex}};}] {
    
     \node (n1) at (0,0)  {$h_{BC}^{0,0}$} ;
     \node (n2) at (1,1) {$h_{BC}^{0,1}$};
     \node (n3) at (-1,1)  {$h_{BC}^{1,0}$};
     \node (n4) at (0,2)  {$h_{BC}^{1,1}$};
     \node (n5) at (2,2)  {$h^{0,2}_{BC}$} ;
     \node (n6) at (-2,2) {$h^{2,0}_{BC}$};
     \node (n7) at (-1,3)  {$h^{2,1}_{BC}$};
     \node (n8) at (1,3)  {$h^{1,2}_{BC}$};
     \node (n8) at (0,4)  {$h^{2,2}_{BC}$};

    \draw [->] (3,2) -- (4,2);

    \node (n1) at (7,0)  {\pgfuseimage{g_tick}} ;
    \node (n2) at (8,1) {\pgfuseimage{g_tick}};
    \node (n3) at (6,1)  {\pgfuseimage{g_tick}};
    \node (n4) at (7,2)  {\pgfuseimage{g_tick}};
    \node (n5) at (9,2)  {\pgfuseimage{g_tick}} ;
    \node (n6) at (5,2) {\pgfuseimage{g_tick}};
    \node (n7) at (6,3)  {\pgfuseimage{r_cross}};
    \node (n8) at (8,3)  {\pgfuseimage{r_cross}};
    \node (n8) at (7,4)  {\pgfuseimage{g_tick}};
}

\end{figure}
\end{center}

It would certainly be very interesting to know how the Hodge numbers would behave for higher dimensional almost complex manifolds. 

\section{$\boldsymbol{h^{p,q}}$ on the Kodaira-Thurston manifold} \label{KTex}
We summarize the calculation of Hodge numbers for a family of almost complex structures on the Kodaira-Thurston manifold. 
\subsection{Almost complex structures on the Kodaira-Thurston manifold}
The Kodaira-Thurston manifold $\KT$ is defined to be the direct product $S^{1}\times (H_{3}(\mathbb{Z})\backslash H_{3}(\mathbb{R}))$, where $S^1$ is parametrized by $t\in \mathbb R$ and $H_{3}(\mathbb{R})$ denotes the Heisenberg group
$$H_{3}(\mathbb{R}) = \left\{ \begin{pmatrix}
    1 & x & z\\
    0 & 1 & y\\
    0 & 0 & 1
\end{pmatrix}\in GL(3,\mathbb{R}) \right\} $$
and $H_{3}(\mathbb{Z})$ is the subgroup  $H_3(\mathbb R)\cap GL(3,\mathbb{Z})$ acting on $H_{3}(\mathbb{R})$ by left multiplication.
\bigskip

Vector fields $\der{t}, \der{x}, \der{y}+x\der{z}$ and $\der{z}$ are well defined, and form a basis at each point. With respect to this basis, a $4\times 4$ matrix whose square is $-Id$ would define an almost complex structure on $\KT$. In this section, we are mainly interested in  a family of non-integrable almost complex structures defined by 
$$J_{a,b} = \begin{pmatrix}
        0 & -1 & 0 &  0\\
        1 &  0 & 0 &  0\\
        0 &  0 & a &  b\\
        0 &  0 & c &  -a
        \end{pmatrix}$$
with $c = -\frac{a^{2}+1}{b}, a, b\ne 0 \in \mathbb R$.
$$V_{1} = \frac 12 \left(\der{t}-i\der{x}\right)\quad \mathrm{\&}\quad V_{2} = \frac 12 \left(\left(\der{y}+x\der{z}\right)-\frac{a-i}{b}\der{z}\right) $$
span $T^{1,0}_{x}M$ at every point, along with their dual 1-forms
$$\phi_{1} = \dif{t}+i\dif{x}\quad\mathrm{\&}\quad \phi_{2} = (1-ai)\dif{y}-ib(\dif{z}-x\dif{y}). $$
We can see the non-integrability of $J_{a,b}$ by observing that $$\dif{\phi_2}=ib \dif{x}\wedge \dif{y}=\frac{b}{4}(\phi_1\wedge\phi_2+\phi_1\wedge\bar{\phi}_2-\bar{\phi}_1\wedge \phi_2-\bar{\phi}_1\wedge\bar{\phi}_2)$$ has a nontrivial $(0,2)$ part.

\subsection{Hodge numbers  $\boldsymbol{h^{1,1}}$, $\boldsymbol{h^{1,0}}$ and $\boldsymbol{h^{2,0}}$} \label{111020}
By Serre duality, we only need to determine the Hodge numbers $h^{1,0}$, $h^{2, 0}$, $h^{0,1}$, and $h^{1,1}$. 

 First, every $J_{a,b}$ is almost K\"ahler with compatible symplectic structure $\omega_{a,b}=\frac{i}{2}(\phi_1\wedge\bar\phi_1+\phi_2\wedge\bar\phi_2)=dt\wedge dx+bdz\wedge dy$. We call the almost K\"ahler metric determined by $J_{a,b}$ and $\omega_{a,b}$ a ``{\it standard orthonormal metric}". 
 For $J_{a,b}$ with any compatible almost K\"ahler metric, we have $h^{1,1}=3$ by Theorem \ref{akh11}. It implies $\tilde h^{1,1}=3$  for $J_{a,b}$.

We have seen that $h^{2,0}$ and $h^{1,0}$ are almost complex invariants. 
For $J_{a,b}$, $h^{1,0}$ and $h^{2,0}$ are calculated in \cite{CZ}
$$h^{2,0} =
\begin{cases}
    1 & b \in 4\pi \mathbb{Z}, b\neq 0\\
    0 & b \not\in 4\pi \mathbb{Z}
\end{cases},$$
$$h^{1,0}=1.$$

The method of computation is the classical Fourier series. We compute $h^{2,0}$ in the following. The computation of $h^{1,0}$ is similar. 

Write $s\in \Omega^{2,0}(M)$ as $f\phi_1\wedge\phi_2$. So $$\bar\partial s=(\bar V_1(f)+f\frac{b}{4})\phi_1\wedge\bar\phi_1\wedge\phi_2+\bar V_2(f)\phi_1\wedge \phi_2\wedge\bar\phi_2.$$

As $\bar\p^* s=0$ automatically, $s$ is a $\bar\p$-harmonic $(2,0)$-form if and only if
\begin{align}\label{04}\bar V_1(f)+\frac{b}{4}f&=0,\\
\label{5}\bar V_2(f)&=0.
\end{align}

Apply $V_2$ to \eqref{5} and write $f=u+iv$ we have 

\begin{align}\label{8}
 &\frac{\p^2 u}{\p y^2}+2Re(\lambda)\frac{\p^2u}{\p y\p z}+|\lambda|^2 \frac{\p^2u}{\p z^2}=0,\\
\label{9}
&\frac{\p^2 v}{\p y^2}+2Re(\lambda)\frac{\p^2v}{\p y\p z}+|\lambda|^2 \frac{\p^2v}{\p z^2}=0,
\end{align}
where $\lambda=x+\frac{c}{a+i}$.

The $\KT$ has a structure of $T^2$-bundle over $T^2$ by $$\pi(t,x, y, z)=(t,x).$$
Restricting $f$ to a single fibre \eqref{8} and \eqref{9} become elliptic since $Re(\lambda)^2\le |\lambda|^2$. As the fibre is compact, by the maximal principle, $f$ is constant with respect to $y$ and $z$. 

Viewing $f$ as a function on the base $T^2$, we have the Fourier series $$\mathcal{F}(f)=\sum_{(k,l)\in \mathbb{Z}^2}f_{k,l}e^{2\pi i(kt+lx)},\ f_{k,l}=\int_{T^2} f(t,x)e^{-2\pi i(kt+lx)}\dif{t}\dif{x}.$$
 Equation \eqref{04} becomes
$$\sum_{(k,l)\in \mathbb{Z}^2} (\frac{b}{4}+\pi(ik-l))f_{k,l} e^{2\pi i(kt+lx)}=0.$$
If $b\notin 4\pi\mathbb{Z}$, then $\frac{b}{4}+\pi(ik-l)\neq 0$ for any $(k,l)\in \mathbb{Z}^2$. So $f_{k,l}=0$ and $f=0$. If $b=4l\pi$ for some $l\in\mathbb{Z}\backslash \{0\}$, then $f=Ce^{2\pi i lx}$ are the solutions. Therefore we get
$$h^{2,0}=\begin{cases} 0, b\notin 4\pi\mathbb{Z}\\
1, b\in 4\pi\mathbb{Z}\end{cases}$$

\subsection{Hodge number $\boldsymbol{h^{0,1}}$}
The method of computing $h^{0,1}$ is very different from above. First, we need to specify an almost Hermitian metric. We choose the standard orthonormal metric which is defined in Section \ref{111020}. 

Write a $(0,1)$-form $s = f\bar{\phi}_{1}+g\bar{\phi}_{2}$ with $f,g\in C^{\infty}(\KT)$. $s$ is $\bar\p$-harmonic if and only if $\bar{\partial}s = 0$ and $\partial *s = 0$, which is equivalent to \begin{align}
    -\bar{V}_{2}(f)+\bar{V}_{1}(g)+g\frac b4&=0,\label{V1}\\
    V_{1}(f)+V_{2}(g)&=0.\label{V2}
\end{align}
We are not able to use the classical Fourier series as in the computation of $h^{2,0}$, as $f$ and $g$ will depend on all the variables in general and the Heisenberg group is non-abelian.

However, theoretically every locally compact group has a Fourier theory since its essence is to decompose Hilbert function spaces with respect to irreducible unitary representations, in the spirit of Peter-Weyl. 
In our situation, the Kodaira-Thurston manifold is derived from the Heisenberg group, whose irreducible unitary representations are classified by the classical Stone-von Neumann theorem.

\section{Solving equations using non-abelian Fourier theory}\label{pdeodent}
To solve Equations \eqref{V1} and \eqref{V2}, our idea is to decompose the space of smooth functions (or $L^2$-functions) into smaller spaces of smooth (resp. $L^2$) functions which are preserved under left invariant differential operators $V_1, V_2$, then solve the equations at these smaller spaces. This is essentially a process of non-abelian Fourier theory. This method was used to study geometric problems on nilmanifolds in {\it e.g.} \cite{AT, DS, Ric} and recently \cite{RS, AT21}. The new ingredients in \cite{HZ, HZ2} are the roles played by Stokes phenomenon and Gauss circle problem. 

\subsection{Decomposition of functions}
We talk about a slightly more general setting than just $\KT$ following \cite{Ho}. Let $M$ be a $K$-dimensional torus bundle over circle, which is described as a mapping torus determined by $A\in GL(K, \mathbb Z)$. The lattice $\mathbb Z^K$ could be partitioned into orbits $Orb_{\bf y}=\{(A^T)^{\xi}\bf y| \xi\in \mathbb Z\}$ of action $A^T$. The set of all such orbits are denoted by $\mathcal O$. 

Then any smooth function $f$ on $M$ can be written as an absolutely convergent sum of smooth functions
$$F_{\bf y}(x, {\bf x})=\sum_{\xi\in \mathbb Z}f_{\bf y}(x+\xi) e^{2\pi i{\bf y}\cdot A^{\xi}{\bf x}}, \quad G_{l, {\bf y}}(x, {\bf x})=f_{l, {\bf y}}e^{2\pi i\frac{lx}{N}}\sum_{\xi=1}^Ne^{2\pi i(\frac{l\xi}{N}+{\bf y}\cdot A^{\xi}{\bf x})}$$ as follows
\begin{equation}\label{smoothdec}
f(x, {\bf x})=\sum_{Orb_{\bf y}\in \mathcal O, |Orb_{\bf y}|=\infty} F_{\bf y}(x, {\bf x})+\sum_{Orb_{\bf y}\in \mathcal O, |Orb_{\bf y}|=N<\infty, l\in \mathbb Z}G_{l, {\bf y}}(x, {\bf x}) 
\end{equation}
where $f_{\bf y}(x)\in \mathcal S(\mathbb R)$ (space of Schwartz functions) and $f_{l, {\bf y}}\in \mathbb C$ are given by 

$$f_{\bf y}(x)=\int_{[0,1]^K}f(x, {\bf x})e^{-2\pi i{\bf y}\cdot {\bf x}}d{\bf x}, \quad f_{l, {\bf y}}=\frac{1}{N}\int_0^N f_{\bf y}(x)e^{-\frac{2\pi ilx}{N}}dx.$$

In the case of $\KT$ we have $K=3$ and $A=\begin{pmatrix}
    1 & 0 & 0\\
    0 & 1 & 1\\
    0 & 0 & 1
\end{pmatrix}$. If we write ${\bf x}=(t, y, z)$ and ${\bf y}=(k, m, n)$, we know $|Orb_{\bf y}|<\infty$ if and only if $n=0$. These finite orbits have $N=1$ and are parametrized by $k, l, m\in \mathbb Z$. 
The infinite orbits are parametrized by $k, m, n\in \mathbb Z$ with $0\le m<|n|$. Then the $L^2$ completion of \eqref{smoothdec} is the decomposition 

\begin{equation} L^2(\KT)=\left(\widehat\bigoplus_{n\in \mathbb{Z}\backslash \{0\}}\widehat\bigoplus_{\substack{k,m\in \mathbb{Z}\\0\leq m <|n|}}\mathcal{H}^{k,m,n}\right) \oplus\left(\widehat\bigoplus_{k,l,m\in \mathbb{Z}}\mathcal{H}^{k,l,m,0}\right) \label{decomp}\end{equation}
where
 \begin{equation*}\mathcal H^{k, m, n}=\{ F_{k,m,n} = \sum_{\xi\in \mathbb{Z}}f_{k,m,n}(x+\xi)e^{2\pi i(kt+(m+n\xi)y+nz)} \,|\,f_{k,m,n}\in L^2(\mathbb{R})\}, \end{equation*}

$$\mathcal{H}^{k,l,m,0} = \{F_{k,l,m,0}=f_{k,l,m,0}e^{2\pi i (kt+lx+my)}\,|\,f_{k,l,m,0}\in \mathbb{C}\}. $$

The map $W_{k,m,n}:f_{k, m, n}\mapsto F_{k,m,n}  $ is called the Weil-Brezin map. It transforms $L^2(\mathbb R)$ onto $\mathcal H^{k, m, n}$ and Schwartz functions $\mathcal S(\mathbb R)$ to  $\mathcal H_{k, m, n}$.

From representation theory perspective, $L^2(\KT)$ corresponds to ind$_{\Gamma}^G(1)$ where $\Gamma=\mathbb Z\times H_3(\mathbb Z)$ and $G=\mathbb R\times H_3(\mathbb R)$, and each $\mathcal H^{k, m, n}$ corresponds to the representation $\rho_n$ of the Heisenberg group. 
Moreover, each $\mathcal{H}^{k,l,m,0}$ corresponds to the irreducible unitary representation $\sigma_{lm}$ of the Heisenberg group in classical notation ({\it c.f.} \cite{Fol, stein}). 
In other words, all the direct summands in the decomposition are exactly all the irreducible unitary representations of the Heisenberg group invariant under $\Gamma$ by virtue of the Stone-von Neumann theorem (It is like how $e^{inx}$ sit in $e^{itx}$ in the classical 1D Fourier series-transform situation). This interpretation should also work for other choice of matrix $A$ in the more general decomposition \eqref{smoothdec}.

\subsection{Solving PDEs using decomposition I: ODEs for $\mathcal H_{k,m,n}$}\label{nne0}

Since left-invariant differential operators preserve each $\mathcal H_{k,m,n}= \mathcal H^{k, m, n}\cap C^{\infty}(\KT)$,  we can solve the equations for $\bar\p$-harmonic $(0,1)$-forms \eqref{V1} and \eqref{V2} on each $\mathcal H_{k,m,n}$ or $\mathcal H^{k,l,m,0}$ and take linear combinations. The treatment for $n\ne 0$ and $n=0$ are different. 

Solutions in $\mathcal{H}_{k,m,n}$, with fixed $n\neq0$ and $0\leq m<|n|$, take the form of
\begin{align*}
    F_{k,m,n} = \sum_{\xi\in \mathbb{Z}}f_{k,m,n}(x+\xi)e^{2\pi i(kt+(m+n\xi)y+nz)},\\
    G_{k,m,n} = \sum_{\xi\in \mathbb{Z}}g_{k,m,n}(x+\xi)e^{2\pi i(kt+(m+n\xi)y+nz)}.
\end{align*}
Plugging these into the original PDEs we obtain a system of first order ODEs on $\mathbb{R}$.
\begin{align}
    \frac{d}{dx}
    \begin{pmatrix}
        f_{k,m,n}\\
        g_{k,m,n}
    \end{pmatrix}
    = (A_{n}x+B_{k,m,n})
    \begin{pmatrix}
        f_{k,m,n}\\
        g_{k,m,n}
    \end{pmatrix}\label{Main DE}
\end{align}
with
$$ A_{n} = 2\pi\begin{pmatrix}
        0 & n\\
        n & 0\\
        \end{pmatrix},$$
$$ B_{k,m,n} = 2\pi\begin{pmatrix}
        k                   & m-n\frac{a-i}{b}\\
        m-n\frac{a+i}{b}    & \frac{b}{4\pi}i-k\\
        \end{pmatrix}.$$

We will have two independent pairs of solutions to each of the above ODE \eqref{Main DE}, however these solutions may not give rise to smooth functions on $\KT$ through Weil-Brezin transform. In particular, any solution with $f_{k,m,n}(x)$ or $g_{k,m,n}(x)\notin L^{2}(\mathbb{R})$ will not produce a valid $F_{k,m,n}$ and $G_{k,m,n}$ as the series will not converge.

By basic ODE theory, as $x\rightarrow +\infty$  we have two independent local solutions, one that grows like $e^{|n|\pi x^{2}}$ and one that decays like $e^{-|n|\pi x^{2}}$, and likewise for large negative $x$. If we have a single solution that decays in both directions then it must be Schwartzian. However, in most cases, a solution decays in one direction may grow in the other direction.  

The study of these end behaviours is eventually the ``Stokes phenomenon". The following is sufficient for our use. 

\begin{thm}[Holt-Z. \cite{HZ}]\label{1stode}
Let $A,B \in M_{2}(\mathbb{C})$ be matrices and let $A$ have two distinct, real eigenvalues $\lambda_{1}$, $\lambda_{2}$ with $\lambda_{1}>0>\lambda_2$ then the equation
\begin{align}
    \frac{d}{dx}
    \begin{pmatrix}
        f\\
        g
    \end{pmatrix}
    = (Ax+B)
    \begin{pmatrix}
        f\\
        g
    \end{pmatrix}\label{Ax+B}
\end{align}
has a pair of solutions $f,g\in L^{2}(\mathbb{R})$ if and only if the following holds: Given $T\in GL(2, \mathbb C)$ such that $TAT^{-1}$ is diagonal and writing $TBT^{-1}$ as $\begin{pmatrix}
    b_{1} & b_{2}\\
    b_{3} & b_{4}
\end{pmatrix}$ we have $b_{2}b_{3}\in (\lambda_{1}-\lambda_{2})\cdot \mathbb{Z}^-$, and in this situation both $f$ and $g$ are Schwartz functions.
\end{thm}

Without loss, assume $n>0$. We have $\lambda_{1} = 2\pi n, \lambda_{2} = -2\pi n$ and
$$TB_{k, m, n}T^{-1} = 2\pi\begin{pmatrix}
      m-\frac{na}{b} +\frac{b}{8\pi}i & k-\frac nb i-\frac{b}{8\pi}i\\
    k+\frac{n}{b}i-\frac{b}{8\pi}i & -m+\frac{na}{b}+\frac{b}{8\pi}i
\end{pmatrix}, $$
where $T =\frac{1}{\sqrt{2}}\begin{pmatrix}
    1 & 1\\
    1 & -1
\end{pmatrix}$. 
In order for us to have a pair of solutions $f_{k,m,n}, g_{k,m,n} \in L^{2}(\mathbb{R})$ we must have 
$$\frac{\pi^2}{n}\left(k-\frac nb i-\frac{b}{8\pi}i\right)\left(k+\frac{n}{b}i-\frac{b}{8\pi}i\right)\in 4\pi \mathbb{Z}^-. $$

This holds only when $8\pi d^2\in \mathbb Z[\sqrt{D}]$ for some integer $D>0$. Here $d=\frac{b}{8\pi}$. In particular, no rational number $d=\frac{p}{q}$ would satisfy it for any choice of $u, n\in \mathbb Z$.

\subsection{Solving PDEs using decomposition II: lattice points counting for $\mathcal H_{k,l,m,0}$} \label{n=0lattice}

The solutions  in $\mathcal{H}_{k,l,m,0}$ take the form
\begin{align*}
    F_{k,l,m,0} = f_{k,l,m,0}e^{2\pi i (kt+lx+my)},\\
     G_{k,l,m,0} = g_{k,l,m,0}e^{2\pi i (kt+lx+my)}.
\end{align*}

We have a solution for \eqref{V1} and \eqref{V2}  when $f_{k,l,m,0}$ and $g_{k,l,m,0}$ satisfy
\begin{align*}
    -mf_{k,l,m,0}+\left(k+il-\frac{bi}{4\pi}\right)g_{k,l,m,0}&=0,\\
    (k-il)f_{k,l,m,0}+mg_{k,l,m,0}&=0.
\end{align*}

When $m=0$, 
we have solutions 
\begin{equation}f=C_{0},\quad g=0, \label{sol 1}\end{equation}
\begin{equation} f=0,\quad g=C_{1}e^{2\pi i lx}, \quad \hbox{when } b = 4\pi l \in 4\pi\mathbb{Z}\backslash\{0\}. \label{sol 2}\end{equation}

When $m\ne 0$,   we can rewrite our equations as 
$$ \left(k^{2}+l^{2}+m^{2}-\frac{b}{4\pi}l-\frac{b}{4\pi}ik\right)f_{k,l,m,0}=0,$$
$$ g_{k,l,m,0}=-\frac{k-il}{m} f_{k,l,m,0}.$$
A nontrivial solution exists if and only if
$$ k^{2}+l^{2}+m^{2}-\frac{b}{4\pi}l-\frac{b}{4\pi}ik=0.$$
This is the case exactly when $k=0$ and nonzero $l,m$ are chosen such that $b=4\pi (l^{2}+m^{2})/l$. This yields the solutions
\begin{equation}f=mC_{2}e^{2\pi i(lx+my)},\quad g = ilC_{2}e^{2\pi i(lx+my)}. \label{sol 3}\end{equation}

 The solutions \eqref{sol 2} and \eqref{sol 3} could be combined by taking $m\in \mathbb Z$ in \eqref{sol 3}. 

Now the counting of solutions in $\mathcal H_0$ is equivalent to the number theoretic question (let $d=\frac {b}{8\pi}$):

\bigskip

{\it
How many pairs of integers $(m,l)$ satisfy
\begin{equation}\label{circlat}
(l-d)^{2}+m^{2}=d^{2}?
\end{equation}
}


Counting the number of solutions can be thought of as asking how many lattice points lie on a circle with centre $(d,0)$ and radius $d$.
Apparently, \eqref{circlat} has non-trivial solutions only when $d$ is rational (which leads to no solution in $\mathcal H_{k,m,n}$).

When $d$ is an integer this problem is very well understood and the number of such integer pairs is denoted $r_2(d^2)$.
 Write $d^{2}$ as a unique product of prime numbers 
$$ d^{2} = 2^{\alpha_{0}}p_1^{\alpha_{1}}\dots p_s^{\alpha_{s}}q_1^{\beta_{1}}\dots q_t^{\beta_{t}},$$
where $p_{i}\equiv 3 \mod 4$ for all $i$ and $q_{j}\equiv 1 \mod 4$ for all $j$. The number of solutions is then given by 
$$h^{0,1}= r_2(d^2)= 4(\beta_{1}+1)(\beta_{2}+1)\dots (\beta_{t}+1). $$

The Hodge number $h^{0,1}$ of $J_{a, b}$ in this case is equal to the number of lattice points lying on a circle with centre $(d,0)$ and radius $d$.

\begin{thm}[Holt-Z. \cite{HZ}]\label{h01}
When $q\le 5$, we have 

$$
h^{0,1}=\begin{cases} \begin{array}{cll}
4(\beta_{1}+1)(\beta_{2}+1)\dots (\beta_{t}+1)& \hbox{ if $q=1$},\\
2(\beta_{1}+1)(\beta_{2}+1)\dots (\beta_{t}+1)& \hbox{ if $q=2$},\\
(\beta_{1}+1)(\beta_{2}+1)\dots (\beta_{t}+1)& \hbox{ if $q=3$},\\
(\beta_{1}+1)(\beta_{2}+1)\dots (\beta_{t}+1)& \hbox{ if $q=4$},\\
(\beta_{1}+1)(\beta_{2}+1)\dots (\beta_{t}+1)& \hbox{ if $q=5$}.
\end{array}
\end{cases}
$$
when $\frac{p}{q}=d=\frac{b}{8\pi}\ne 0$, $\gcd(p,q)=1$, and $p^2= 2^{\alpha_{0}}p_1^{\alpha_{1}}\dots p_s^{\alpha_{s}}q_1^{\beta_{1}}\dots q_t^{\beta_{t}}$ where $p_{i}\equiv 3 \mod 4$ for all $i$ and $q_{j}\equiv 1 \mod 4$ for all $j$.
\end{thm}

Combination of the discussion of Sections \ref{nne0} and \ref{n=0lattice} leads to Theorems \ref{intro1}. 

Our method to find all the $\bar\partial$-harmonic forms could be summarised by the following diagram.

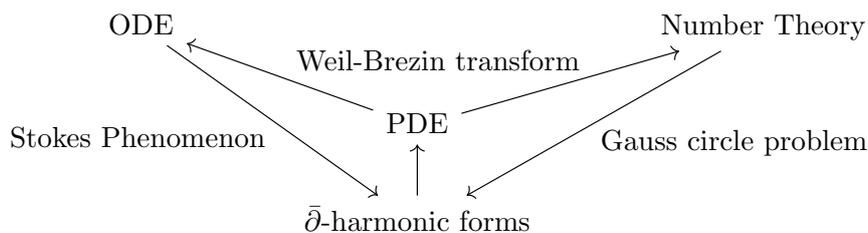
\begin{center}
\begin{tikzcd}[column sep=large]
\hbox{ODE} \arrow{ddr}[swap]{\hbox{Stokes Phenomenon}}
&   &\hbox{Number Theory}\arrow{ddl}{\hbox{Gauss circle problem}}\\
                                       &\hbox{PDE}\arrow{ur} \arrow{ul}[swap]{\hbox{ Weil-Brezin transform}}&\\
                                      &\bar\partial\hbox{-harmonic forms}\arrow{u}&
\end{tikzcd}
\end{center}

\subsection{Kodaira-Spencer question}
We introduce an almost K\"ahler solution of Question \ref{KS} following \cite{HZ2}. A non almost K\"ahler solution is described in \cite{HZ} by a similar calculation.

We consider a family of Hermitian metrics $h_{a,b,\rho}$ compatible with $J_{a,b}$, where we take $\rho \in \mathbb{R}^+$,
$$h_{a,b,\rho} = 2(\phi^1 \otimes \bar{\phi}^1 + \rho \phi^2 \otimes \bar{\phi}^2 ). $$
The associated fundamental form $$\omega= -\frac{i}{2}( h - \bar h) =-2i(\phi^1 \wedge \bar{\phi}^1 + \rho\phi^2 \wedge \bar{\phi}^2 )=4(dx \wedge dt+\rho dy\wedge (dz-xdy))$$ is independent of $a,b$ and closed. Hence, the metrics $h_{a,b,\rho}$ are almost K\"ahler.  

Then the equations for a $\bar\partial$-harmonic form $s=f\bar \phi_1+g\bar\phi_2$ are
\begin{align}
-\bar{V}_2 (f) + \bar{V}_1 (g) +g \frac b4 &=0, \label{V1+}\\
\rho V_1 (f) +V_2 (g) &=0. \label{V2+}
\end{align}

Apply the decomposition and first look for solutions for a fixed $n\neq0$, we obtain
$$\frac{d}{dx} \begin{pmatrix}
    f_{k,m,n}\\
    g_{k,m,n}
\end{pmatrix}=
(Ax+B)\begin{pmatrix}
    f_{k,m,n}\\
    g_{k,m,n}
\end{pmatrix}$$
where
$$A = 2\pi n\begin{pmatrix}
    0   &   \frac{1}{\rho}\\
    1   &   0
\end{pmatrix}, \quad \quad 
B = 2\pi\begin{pmatrix}
    k   &   \frac {1}{\rho}\left(m-\frac{a-i}{b}n\right)\\
    m - \frac{a+i}{b}n    &   i\frac{b}{4\pi}-k
\end{pmatrix}.$$

Apply Theorem \ref{1stode}, we know that there is no solution unless $8\pi \sqrt \rho d^2$ is a quadratic integer. In particular, $\sqrt \rho d^2 \in \frac{1}{\pi}\bar{\mathbb{Q}}$, where $\bar{\mathbb{Q}}$ is used to denote the algebraic numbers.

The equations for $n=0$ would be
$$-m f_{k,l,m,0} + \left(k+il -i\frac{b}{4\pi} \right)g_{k,m,n}=0 $$
$$\rho \left(k-il \right)f_{k,m,n} + mg_{k,m,n} =0$$
Setting $m=0$ would lead to the solutions
$$f=C_{1}, g=0 \quad \text{and} \quad f=0, g=C_{2}e^{2\pi i \frac{b}{4\pi}x} $$
for any $C_{1}, C_{2} \in \mathbb{C}$. 

When $m\ne 0$, there will be solutions only when $k=0$ and non-zero $l,m\in \mathbb{Z}$ are chosen such that
\begin{equation}\label{gauss'}
\left(\frac{m}{\sqrt{\rho}}\right)^2 +(l-d)^2 = d^2. 
\end{equation}
Such a choice for $l,m$ will give us the solutions
$$f= mC_{3}e^{2\pi i (lx+my)}, g=i\rho l C_{3}e^{2\pi i (lx+my)}  $$
$C_{3}\in\mathbb{C}$. 
Equation \eqref{gauss'} is equivalent to asking how many points on the lattice given by $\mathbb{Z}\times\frac{1}{\sqrt \rho}\mathbb{Z}$ intersect a circle of radius $d$ and centre $(d,0)$. Note that the two solutions found when $m=0$ can be considered to correspond to the lattice points $(0,0)$ and $(\frac{b}{4\pi},0)$. 

Consider the case when $\frac{b}{8\pi}=d=1$. Bringing together the above discussion we find that

$$ h^{0,1} = \begin{cases} 
      4 & \sqrt \rho \in \mathbb{Z} ,\\
      2 & \sqrt \rho \in \mathbb{Q}\backslash \mathbb{Z} .
   \end{cases} $$
   That is, $h^{0,1}$ varies under almost K\"ahler metric deformations.

\subsection{Solving Stokes phenomenon}
The ODE Theorem \ref{1stode} is crucial for this ``PDE to ODE and NT" method. We discuss its proof in this section. 

The trick is to make a substitution $$    \begin{pmatrix}
        \psi\\
        \phi
    \end{pmatrix}=  e^{-\frac{1}{2}\lambda_{2}x^{2}}T
    \begin{pmatrix}
        f\\
        g
    \end{pmatrix},$$

 such that $A$ becomes a matrix with only one non-zero entry 
\begin{align*} 
    \frac{d}{dx}
    \begin{pmatrix}
        \psi\\
        \phi
    \end{pmatrix}
    = \left(\begin{pmatrix}
        \lambda_{1}-\lambda_{2} & 0\\
        0 & 0
    \end{pmatrix}x+TBT^{-1}\right)
    \begin{pmatrix}
        \psi\\
        \phi
    \end{pmatrix}.
\end{align*}

Then both $\phi$ and $\psi$ will satisfy (different) second order ODEs of type
$$(p_{2}x+q_{2})h''+(p_{1}x+q_{1})h'+(p_{0}x+q_{0})h=0, $$
which could be solved using a Laplace integral transform.

We can write $h$ as
$$h(x) = \int_{C}\varphi(s) e^{sx}ds $$
where $C$ is some contour in the complex $s$-plane. Then, defining
\begin{align*}
    P(s) = p_{2}s^{2}+p_{1}s+p_{0}\\
    Q(s) = q_{2}s^{2}+q_{1}s+q_{0}
\end{align*}
and choosing $C$ so that
$$V(s) = \exp\left(\int^{s} \frac{Q(\sigma)}{P(\sigma)}d\sigma\right)e^{s x} $$
takes the same value at both (possibly infinite) endpoints for all $x$ when $s$ parameterises the contour $C$, we can find a solution
$$\varphi(s) = \frac{1}{P(s)}\exp\left(\int^{s} \frac{Q(\sigma)}{P(\sigma)}d\sigma\right) \quad \mathrm{\&}\quad h(x)=\int_C\frac{V(s)}{P(s)}ds. $$

In our specific case, $$P_{\phi}(s) = (\lambda_{1}-\lambda_{2})(b_{4}-s), \quad \quad Q_{\phi}(s) = s^{2}-(b_{1}+b_{4})s+(b_{1}b_{4}-b_{2}b_{3}),$$
which gives us the solution
\begin{align*}
    \phi(x) = \frac{1}{\lambda_2-\lambda_1}\int_{C}(s-b_{4})^{\frac{b_{2}b_{3}}{\lambda_{1}-\lambda_{2}}-1}\exp(-\frac{1}{\lambda_{1}-\lambda_{2}}\left(\frac{s^2}{2}-b_{1}s\right)+xs) ds\label{phiint}
\end{align*}
with 
$$V_{\phi}(s) = (s-b_{4})^{\frac{b_{2}b_{3}}{\lambda_{1}-\lambda_{2}}}\exp(-\frac{1}{\lambda_{1}-\lambda_{2}}\left(\frac{s^2}{2}-b_{1}s\right)+xs). $$

$V_{\phi}(s)$ tends to zero as $s$ grows large within the shaded regions. 

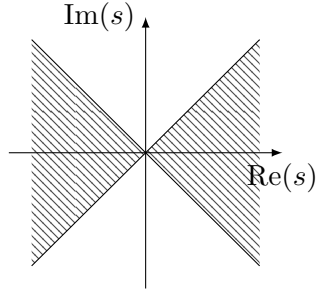
\begin{center}
\begin{tikzpicture}
\fill[pattern color=black!70,scale=.3,domain=.2:5.5,smooth,pattern=north west lines] (-5,5) -- plot({\x},-{\x}) -- (5,-5) -- (5,5) -- plot ({\x},{\x}) -- (-5,-5) -- cycle;
    \draw[scale=.3,domain=-5:5,smooth] plot ({\x},{\x});
    \draw[scale=.3,domain=-5:5,smooth] plot ({\x},{-\x});
    \draw[>=latex,->] (-1.8,0) -- (1.8,0) node[below] {$\Re(s)$};
    \draw[>=latex,->] (0,-1.8) -- (0,1.8) node[left] {$\Im(s)$};
\end{tikzpicture}
\end{center}

So we can choose the contours

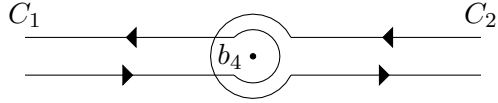
\begin{center}
\begin{tikzpicture}    
\draw[fill] (0,0) circle (1pt) node[left] {$b_{4}$};
\begin{scope}[every node/.style={sloped,allow upside down}]
    \draw (-.25,.25) -- node {\midarrow} (-3,.25) node[above] {$C_{1}$};
    \draw (-3,-.25) -- node {\midarrow} (-.25,-.25);
    \draw (3,.25)  node[above] {$C_{2}$} -- node {\midarrow} (.5,.25);
    \draw (.5,-.25) -- node {\midarrow} (3,-.25);
\end{scope}
    \draw [domain=-135:135] plot ({0.3535*cos(\x)}, {0.3535*sin(\x)}); 
    \draw [domain=26.45:333.55] plot ({0.559*cos(\x)}, {0.559*sin(\x)});
\end{tikzpicture}
\end{center}

In fact, we were looking at contours $\tilde{C}_{1}$ and $\tilde{C}_{2}$ that is transformed from $C_1$ and $C_2$ by substitution $t := s-(b_{1}+(\lambda_{1}-\lambda_{2})x)$.
When $\frac{b_{2}b_{3}}{\lambda_{1}-\lambda_{2}}\notin\mathbb{Z}$, $C_1$ and $C_2$ give two independent solutions $\phi_{1}$ and $\phi_{2}$.
Analyzing this end behaviors, we know

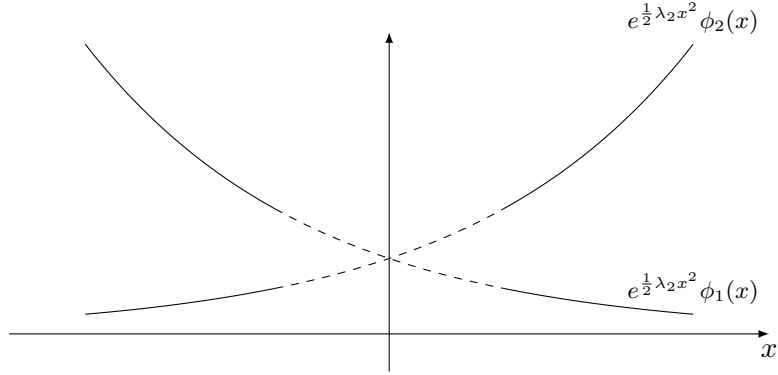
\begin{center}
    \begin{tikzpicture}
        \draw [domain=-4:-1.5] plot ({\x},{1.4^(\x)}) ;
        \draw [domain=-4:-1.5] plot ({\x},{1.4^(-\x)}) ; 
        \draw [domain=-1.5:1.5,dashed]plot ({\x},{1.4^(\x)}) ;
        \draw [domain=-1.5:1.5,dashed] plot ({\x},{1.4^(-\x)}) ;    
        \draw [domain=1.5:4]plot ({\x},{1.4^(\x)}) node[above] {\footnotesize{$e^{\frac{1}{2}\lambda_{2}x^{2}}\phi_{2}(x)$}};
        \draw [domain=1.5:4] plot ({\x},{1.4^(-\x)}) node[above] {\footnotesize{$e^{\frac{1}{2}\lambda_{2}x^{2}}\phi_{1}(x)$}};
        \draw[>=latex,->] (-5,0) -- (5,0) node[below] {$x$};
        \draw[>=latex,->] (0,-.5) -- (0,4);
    \end{tikzpicture}
\end{center}

In particular, they are not $L^2$.

If $\frac{b_{2}b_{3}}{\lambda_{1}-\lambda_{2}}\in\mathbb{Z}^-\cup\{0\}$, then the integrals along the horizontal directions of the path of integration cancel, and the two integrals along $\tilde C_1$ and $\tilde C_2$ reduce to an integral  along a loop around $t_0=b_4-b_1-(\lambda_1-\lambda_2)x$. 
This gives rise to a solution $\phi(x)$ which grows at most as $e^{Kx}$ at both ends. Hence $e^{\frac{1}{2}\lambda_{2}x^{2}}\phi$ decays as $e^{\frac{1}{2}\lambda_{2}x^{2}}$ at both ends.

The situation of $\frac{b_{2}b_{3}}{\lambda_{1}-\lambda_{2}}\in\mathbb{Z}^+$ is slightly more delicate. The reader can check \cite{HZ} for the proof. 

\subsection{Applying ``PDE to ODE and NT" method in other scenarios}
We end this article by applying this method in a couple more scenarios.
\subsubsection{Compute $h^{1,1}$}

The method is also applicable to $h^{1,1}$. We take $J_{a,b}$ and the standard orthonormal metrics for example. We write $(1,1)$-form as  $$s=f^{(1,1)}\phi_{1}\wedge\bar{\phi}_{1}+f^{(1,2)}\phi_{1}\wedge\bar{\phi}_{2}+f^{(2,1)}\phi_{2}\wedge\bar{\phi}_{1}+f^{(2,2)}\phi_{2}\wedge\bar{\phi}_{2}.$$

When $n=0$, we will get three dimensions of solutions $$f^{(1,1)}=C_{0},\quad f^{(1,2)}=C_{1},\quad f^{(2,1)}=-C_{1},\quad f^{(2,2)}=C_{2} .$$

But for $n>0$, we  have a countably many $4\times 4$ ODE systems. We know none of them leads to $\bar\p$-harmonic forms as Theorem \ref{akh11} shows that $h^{1,1}=3$ for any almost K\"ahler metric. But how can we prove it directly by solving equations?

In general, we will need to study Stokes phenomenon for $N\times N$ ODE systems. But our practical knowledge of it seems very limited.

\subsubsection{Euclidean manifolds}
An interesting calculation has been done by Tom Holt for almost complex structures on manifolds with Euclidean geometry in \cite{Ho}. More precisely, he applies it to a $3$-torus bundle with the monodromy a matrix with finite order. The particular example he studied has $$A= \begin{pmatrix}
        0 & 1 &  0\\
        0 &  0 & 1 \\
        1 &  0 & 0
   
        \end{pmatrix}.$$

By applying this method, he can reduce the problem to study the ``Stokes phenomenon" for two sequences $\{a_k\}_{k=-\infty}^{\infty}$ and $\{b_k\}_{k=-\infty}^{\infty}$ with recurrence relation 
$$ \begin{pmatrix}
    a_k\\
    b_k
\end{pmatrix}=
\frac{1}{dk+e}(Ak^2+Bk+C)\begin{pmatrix}
    a_{k-1}\\
    b_{k-1}
\end{pmatrix}.$$

Namely, it is equivalent to asking whether they are in the following discrete Schwartz space $$\mathcal S(\mathbb Z)=\{a_k: \mathbb Z\rightarrow \mathbb C| \sup_{k\in \mathbb Z} |k^pa_k|<\infty, \forall p\in \mathbb N\}.$$

It is interesting to notice that when we look at the  generating functions $A(x)=\sum_{k=-\infty}^{\infty}a_kx^k$ and $B(x)=\sum_{k=-\infty}^{\infty}b_kx^k$, the above recurrence relation would turn to be a second order ODE system on $A(x)$ and $B(x)$. Can this equation be derived from an elliptic equation using PDE to ODE and NT method?

Holt has also done some more computations on manifolds with geometries $Nil^4$ and $Sol^3\times \mathbb E$ in \cite{Hothesis}.

\end{document}